\documentclass{amsart}

\usepackage{fancyhdr}
\usepackage{lastpage}
\usepackage{stmaryrd,yhmath}

\newcommand{\bdis}{\begin{displaymath}}
\newcommand{\edis}{\end{displaymath}}
\newcommand{\be}{\begin{equation}}
\newcommand{\ee}{\end{equation}}
\newcommand{\mbb}{\mathbb}
\newcommand{\mcal}{\mathcal}

\newcommand{\vp}{\varphi}

\newcommand{\mT}{\mathring{T}}

\newcommand{\zf}{\zeta\left(\frac{1}{2}+it\right)}

\theoremstyle{definition}

\newtheorem{cor}[]{Corollary}

\theoremstyle{remark}
\newtheorem{remark}[]{Remark}

\newtheorem*{mydef1}{{\bf Theorem}}

\numberwithin{equation}{section}

\begin{document}

\title{Jacob's ladders and some nonlinear integral equations connected with the Poisson-Lobachevsky integral}

\author{Jan Moser}

\address{Department of Mathematical Analysis and Numerical Mathematics, Comenius University, Mlynska Dolina M105, 842 48 Bratislava, SLOVAKIA}

\email{jan.mozer@fmph.uniba.sk}


\begin{abstract}
We obtain some new properties of the signal generated by the Riemann zeta-function in this paper. Namely, we show the connection between
the function $\zf$ and a nonlinear integral equation related to the Poisson-Lobachevsky integral.
\end{abstract}

\maketitle

\section{The result}

\subsection{}

Let us remind the Poisson parametric integral

\be \label{1.1}
\int_0^\pi \ln(1+a^2-2a\cos x){\rm d}x=\left\{\begin{array}{rcl} 0 & , & a\in (0,1), \\ 2\pi\ln a & , & a>1. \end{array} \right.
\ee

It is well-known that Lobachevsky applied his new geometric system to the calculations of the wide set of integrals. One of these integrals
(see \cite{2}, eq. (111)) leads, after a small transformation, right to the integral (\ref{1.1}).

\begin{remark}
The integral (111) of \cite{2} was used by Lobachevsky to answer the unfair criticism of Ostrogradsky.
\end{remark}

\subsection{}

In this paper we obtain some new properties of the signal

\bdis
Z(t)=e^{i\vartheta(t)}\zf
\edis
generated by the Riemann zeta-function, where

\bdis
\vartheta(t)=-\frac t2\ln\pi+\text{Im}\ln\Gamma\left(\frac 14+i\frac t2\right)=\frac t2\ln\frac{t}{2\pi}-\frac t2-\frac \pi 8 +\mcal{O}\left(\frac 1t\right) ,
\edis

namely, the nonlinear integral equation connected with the function $\zf$ and the Poisson-Lobachevsky integral (\ref{1.1}). \\

Next let us remind that

\bdis
\tilde{Z}^2(t)=\frac{{\rm d}\vp_1(t)}{{\rm d}t},\ \vp_1(t)=\frac 12\vp(t) ,
\edis

where

\be \label{1.2}
\tilde{Z}^2(t)=\frac{Z^2(t)}{2\Phi^\prime_\vp[\vp(t)]}=\frac{\left|\zf\right|^2}{\left\{ 1+\mcal{O}\left(\frac{\ln\ln t}{\ln t}\right)\right\}\ln t}
\ee

(see \cite{3}, (3.9); \cite{4} (1.3); \cite{9}, (1.1), (3.1), (3.2)), and $\vp$ is the Jacob's ladder, i.e. a solution of the nonlinear integral equation
(see \cite{3})

\bdis
\int_0^{\mu[x(T)]} Z^2(t)e^{-\frac{2}{x(T)}t}{\rm d}t=\int_0^TZ^2(t){\rm d}t ,
\edis

(there is an infinite set of Jacob's ladders).

\subsection{}

The following theorem holds true.

\begin{mydef1}
Every Jacob's ladder $\vp_1(t)=\frac 12\vp(t)$, where $\vp(t)$ is the exact solution of the nonlinear integral equation
\bdis
\int_0^{\mu[x(T)]} Z^2(t)e^{-\frac{2}{x(T)}t}{\rm d}t=\int_0^TZ^2(t){\rm d}t ,
\edis
is the asymptotic solution of the following nonlinear integral equation
\be \label{1.3}
\frac
{\int_{x^{-1}(T+\tau_a)}^{x^{-1}(T+\pi)}\ln[1+a^2-2a\cos(x(t)-T)]\left|\zf\right|^2}
{-\int_{x^{-1}(T)}^{x^{-1}(T+\tau_a)}\ln[1+a^2-2a\cos(x(t)-T)]\left|\zf\right|^2}=1 ,
\ee
where $x(t)=x(t;a)$, and
\bdis
\tau_a=\arccos \frac a2,\ a\in (0,1),
\edis
i.e. the following asymptotic formula
\be \label{1.4}
\frac
{\int_{\vp_1^{-1}(T+\tau_a)}^{\vp_1^{-1}(T+\pi)}\ln[1+a^2-2a\cos(x(t)-T)]\left|\zf\right|^2}
{-\int_{\vp_1^{-1}(T)}^{\vp_1^{-1}(T+\tau_a)}\ln[1+a^2-2a\cos(x(t)-T)]\left|\zf\right|^2}=1+\mcal{O}\left(\frac{\ln\ln T}{\ln T}\right) ,
\ee
as $T\to\infty$ holds true for every Jacob's ladder and for every fixed $a\in (0,1)$.
\end{mydef1}

\begin{remark}
There are the fixed point methods and other methods of the functional analysis used to study the nonlinear equations. What can be obtained by using these
methods in the case of the nonlinear integral equation (\ref{1.3})?
\end{remark}

This paper is a continuation of the series \cite{3} - \cite{22}.

\section{Corollaries and remarks}

\subsection{}

It is clear that (\ref{1.3}) is followed by:

\begin{cor}
\be \label{2.1}
\int_{x^{-1}(T)}^{x^{-1}(T+\pi)}\ln[1+a^2-2a\cos(x(t)-T)]\left|\zf\right|^2{\rm d}t=0,
\ee
i.e. we have homogenous nonlinear integral equation.
\end{cor}

Next, from (\ref{2.1}), in the case $a=\frac 1b,\ b>1$ we obtain

\begin{cor}
\bdis
\begin{split}
& \int_{x^{-1}(T)}^{x^{-1}(T+\pi)}\ln[1+b^2-2b\cos(x(t)-T)]\left|\zf\right|^2{\rm d}t = \\
& = 2\ln b\int_{x^{-1}(T)}^{x^{-1}(T+\pi)}\left|\zf\right|^2{\rm d}t ,
\end{split}
\edis
(comp. with the integral (\ref{1.1})).
\end{cor}

\begin{remark}
If we apply the usual way in the case (\ref{1.1}), $a\geq 2$, we obtain the following theorem: every Jacob's ladder $\vp_1(t)=\frac 12\vp(t)$, where
$\vp(t)$ is the exact solution of the nonlinear integral equation

\bdis
\int_0^{\mu[x(T)]} Z^2(t)e^{-\frac{2}{x(T)}t}{\rm d}t=\int_0^TZ^2(t){\rm d}t ,
\edis

is the asymptotic solution of the following nonlinear integral equation

\bdis
\begin{split}
& \int_{x^{-1}(T)}^{x^{-1}(T+\pi)}\ln[1+a^2-2a\cos(x(t)-T)]\left|\zf\right|^2{\rm d}t= \\
& = 2\pi\ln a \ln T,\ a\geq 2 ,
\end{split}
\edis

where $\ln[1+a^2-2a\cos(x(t)-T)]\geq 0$, i.e. the following asymptotic formula

\be \label{2.2}
\begin{split}
& \int_{\vp_1^{-1}(T)}^{\vp_1^{-1}(T+\pi)}\ln[1+a^2-2a\cos(\vp_1(t)-T)]\left|\zf\right|^2{\rm d}t\sim \\
& \sim 2\pi \ln a \ln T,\ T\to\infty
\end{split}
\ee

holds true.
\end{remark}

\subsection{}

On Riemann hypothesis the following Littlewood's estimate

\bdis
0<\gamma'-\gamma<\frac{A}{\ln\ln\gamma},\ \gamma\to\infty
\edis

holds true (see \cite{1}), where $\gamma,\gamma'$ stand for consecutive zeroes of the function $\zf$. Since

\bdis
\frac{1}{\gamma'-\gamma}>\frac{1}{A}\ln\ln\gamma>2,\ \gamma>\gamma_0,
\edis

then we obtain from (\ref{2.1})

\begin{cor} On Riemann hypothesis 
\bdis
\begin{split}
& \int_{\vp_1^{-1}(\gamma)}^{\vp_1^{-1}(\gamma+\pi)}
\ln\left[ 1+\frac{1}{(\gamma'-\gamma)^2}-\frac{2}{\gamma'-\gamma}\cos(\vp_1(t)-\gamma)\right]
\left|\zf\right|^2{\rm d}t\sim \\
& \sim -2\pi\ln(\gamma'-\gamma)\ln \gamma,\ \gamma>\max\{ \gamma_0,T_0[\vp_1]\} .
\end{split}
\edis
\end{cor}

Next, we obtain

\begin{cor}
\bdis
\begin{split}
& \int_{\vp_1^{-1}(p)}^{\vp_1^{-1}(p+\pi)}\ln[1+(p'-p)^2-2(p'-p)\cos(\vp_1(t)-p)]\left|\zf\right|^2{\rm d}t\sim \\
& \sim 2\pi\ln(p'-p)\ln p,\ p\geq T_0[\vp_1] ,
\end{split}
\edis
where $p,p';\ p<p'$ denote the consecutive prime numbers.
\end{cor}

\section{Proof of the Theorem}

\subsection{}

Let us remind that the following lemma holds true (see \cite{8}, (2.5); \cite{9}, (3.3)): for every integrable function (in the Lebesgue sense)
$f(x),\ x\in [\vp_1(T),\vp_1(T+U)]$ we have
\be \label{3.1}
\int_T^{T+U}f[\vp_1(t)]\tilde{Z}^2(t){\rm d}t=\int_{\vp_1(T)}^{\vp_1(T+U)}f(x){\rm d}x,\ U\in \left(\left. 0,\frac{T}{\ln T}\right]\right.
\ee
where
\be \label{3.2}
t-\vp_1(t)\sim (1-c)\pi(t) ,
\ee
$c$ is the Euler's constant and $\pi(t)$ is the prime-counting function. In the case $\mT=\vp_1^{-1}(T),\ \widering{T+U}=\vp_1^{-1}(T+U)$ we
obtain from (\ref{3.1})
\be \label{3.3}
\int_{\vp_1^{-1}(T)}^{\vp_1^{-1}(T+U)}f[\vp_1(t)]\tilde{Z}^2(t){\rm d}t=\int_T^{T+U}f(x){\rm d}x .
\ee

\subsection{}

We obtain from (\ref{1.1}), $0<a<1$

\be \label{3.4}
\int_{\tau_a}^\pi\ln(1+a^2-2a\cos\tau){\rm d}\tau=-\int_0^{\tau_a}\ln(1+a^2-2a\cos\tau){\rm d}\tau ,
\ee

where

\be \label{3.5}
\tau_a=\arccos\frac a2 ,
\ee

and

\be \label{3.6}
\begin{split}
& \ln(1+a^2-2a\cos\tau) >0 ,\ \tau\in (\tau_a,\pi) , \\
& \ln(1+a^2-2a\cos\tau) <0 ,\ \tau\in (0,\tau_a) .
\end{split}
\ee

Putting

\bdis
f(t)=\ln[1+a^2-2a\cos(t-T)],\ U=\pi;\ t=\tau+T
\edis

in (\ref{3.3}), we obtain (see (\ref{3.3}), (\ref{3.4}))

\be \label{3.7}
\begin{split}
& \int_{\vp_1^{-1}(T+\tau_a)}^{\vp_1^{-1}(T+\pi)}\ln[1+a^2-2a\cos(\vp_1(t)-T)]\tilde{Z}^2(t){\rm d}t= \\
& = -\int_{\vp_1^{-1}(T)}^{\vp_1^{-1}(T+\tau_a)}\ln[1+a^2-2a\cos(\vp_1(t)-T)]\tilde{Z}^2(t){\rm d}t
\end{split}
\ee

where

\bdis
\int_T^{T+\pi}\ln[1+a^2-2a\cos(t-T)]{\rm d}t=\int_0^\pi\ln(1+a^2-2a\cos\tau){\rm d}\tau .
\edis

\subsection{}

We obtain by using the mean-value theorem for the integrals in (\ref{3.7}) (see (\ref{1.2}), (\ref{3.6}))

\be \label{3.8}
\begin{split}
& \int_{\vp_1^{-1}(T+\tau_a)}^{\vp_1^{-1}(T+\pi)}\ln[1+a^2-2a\cos(\vp_1(t)-T)]\tilde{Z}^2(t){\rm d}t= \\
& =\frac{1}{\left\{ 1+\mcal{O}\left(\frac{\ln\ln t_1}{\ln t_1}\right)\right\}\ln t_1}
\int_{\vp_1^{-1}(T+\tau_a)}^{\vp_1^{-1}(T+\pi)}\ln[1+a^2-2a\cos(\vp_1(t)-T)]\left|\zf\right|^2{\rm d}t , \\
& \int_{\vp_1^{-1}(T)}^{\vp_1^{-1}(T+\tau_a)}\ln[1+a^2-2a\cos(\vp_1(t)-T)]\tilde{Z}^2(t){\rm d}t= \\
& =\frac{1}{\left\{ 1+\mcal{O}\left(\frac{\ln\ln t_2}{\ln t_2}\right)\right\}\ln t_2}
\int_{\vp_1^{-1}(T)}^{\vp_1^{-1}(T+\tau_a)}\ln[1+a^2-2a\cos(\vp_1(t)-T)]\left|\zf\right|^2{\rm d}t ,
\end{split}
\ee

where $t_1,t_2\in (\vp_1^{-1}(T),\vp_1^{-1}(T+\pi))$ and

\be \label{3.9}
t_1=\vp_1^{-1}(T_1),\ t_2=\vp_1^{-1}(T_2),\ T_1,T_2\in (T,T+\pi) .
\ee

\subsection{}

Next, we obtain from (\ref{3.2}) by (\ref{3.9}) ($t_1\to\infty \ \Leftrightarrow \ T\to\infty$)

\be \label{3.10}
t_1-T_1=\mcal{O}\left(\frac{t_1}{\ln t_1}\right) \ \Rightarrow \ 1-\frac{T_1}{t_1}=\mcal{O}\left(\frac{1}{\ln t_1}\right)\to 0,\ T\to\infty ,
\ee

i.e.

\be \label{3.11}
t_1\sim T_1\sim T,\ t_2\sim T_2\sim T ,
\ee

and (see (\ref{3.10}), (\ref{3.11}); $0<T_1-T, T_2-T<\pi$)

\be \label{3.12}
t_1-T=\mcal{O}\left(\frac{T}{\ln T}\right),\ t_2-T=\mcal{O}\left(\frac{T}{\ln T}\right) .
\ee

Now,

\be \label{3.13}
\ln t_1=\ln T+\mcal{O}\left(\frac{t_1-T}{T}\right)=\ln T+\mcal{O}\left(\frac{1}{\ln T}\right) ,
\ee

and similarly,

\be \label{3.14}
\ln t_2=ln T+\mcal{O}\left(\frac{1}{\ln T}\right) .
\ee

Then the formula (\ref{1.3}) follows from (\ref{3.7}) by (\ref{3.8}), (\ref{3.13}) and (\ref{3.14}).

\thanks{I would like to thank Michal Demetrian for helping me with the electronic version of this work.}


\begin{thebibliography}{29}
\bibitem{18}
J. Moser, `Jacob's ladders and the nonlocal interaction of the function $Z^2(t)$ with the function $\tilde{Z}^2(t)$ on the distance $\sim (1-c)\pi(t)$
for the collections of disconnected sets`, (2010), arXiv: 1007.5147.
\end{thebibliography}
\end{document}